\documentclass[12pt, a4paper]{article}
\usepackage[latin1]{inputenc}
\usepackage[english]{babel}
\usepackage{amsmath}
\usepackage[english]{babel}
\usepackage{indentfirst}
\usepackage{amstext}
\usepackage{enumerate}
\usepackage{amsfonts}
\usepackage{textcomp}
\usepackage{amssymb}
\usepackage[dvips]{graphicx}
\usepackage{setspace}
\usepackage[all]{xy}
\usepackage{color}

\usepackage[
margin=1in,
includefoot,
footskip=30pt,
]{geometry} 

\newcommand {\demo}{\hskip -0.6cm{\bf Proof.  }}
\newcommand {\fim}{\hfill{$\square$}\vskip 1pc}
\newcommand {\F}{\mathbb{F}}
\newcommand {\N}{\mathbb{N}}
\newcommand {\Z}{\mathbb{Z}}
\newcommand {\C}{\mathbb{C}}
\newcommand {\GG}{\mathcal{G}}

\newtheorem{teorema}{Theorem}[section]
\newtheorem{lema}[teorema]{Lemma}
\newtheorem{corolario}[teorema]{Corollary}
\newtheorem{definicao}[teorema]{Definition}
\newtheorem{proposicao}[teorema]{Proposition}
\newtheorem{exemplo}[teorema]{Example}
\newtheorem{remark}[teorema]{Remark}

\begin{document}
\onehalfspace

\title{Properties of the gradings on ultragraph algebras via the underlying combinatorics}

\author{Daniel Gon\c{c}alves\footnote{This author is partially supported by Conselho Nacional de Desenvolvimento Cient\'{i}fico e Tecnol\'{o}gico - CNPq and Capes-PrInt.}  \ and Danilo Royer}

\maketitle

\begin{abstract} There are two established gradings on Leavitt path algebras associated with ultragraphs, namely the grading by the integers group and the grading by the free group on the edges. In this paper, we characterize properties of these gradings in terms of the underlying combinatorial properties of the ultragraphs. More precisely, we characterize when the gradings are strong or epsilon-strong. The results regarding the free group on the edges are new also in the context of Leavitt path algebras of graphs. Finally, we also describe the relation between the strongness of the integer grading on an ultragraph Leavitt path algebra and the saturation of the gauge action associated with the corresponding ultragraph C*-algebra.
    
\end{abstract}

\vspace{1.0pc}

\vspace{1.0pc}

\section{Introduction}

The integer grading on Leavitt path algebras plays a key role in their classification program (which parallels Elliott's classification program for C*-algebras, see for example \cite{H1}). In fact, the determination of a full invariant for $\Z$-graded isomorphism of Leavitt path algebras has still to be determined, but important advances have been made since the graded isomorphism conjecture made by Hazrat in \cite{H1}, see for example \cite{AP, ACort, CHG, HH}. The grading of a Leavitt path algebra is also fundamental in the understanding of its algebraic structure, in particular its ideal structure, see \cite{AHLS, Ranga} for example.  

Given the exposed above, the importance of understanding the grading of a Leavitt path algebra is clear. In particular, the integer grading of a Leavitt path algebra has been studied in \cite{H3}, where the unital Leavitt path algebras which are strongly graded are completely characterized, in \cite{Clark, strongradedleavitt}, where strongly graded Leavitt path algebras are characterized in terms of Condition~(Y), and in \cite{groupgradleavitt}, where it is shown that the Leavitt path algebra associated to a finite graph is epsilon-strongly $\Z$-graded. The grading over the free group of the edges has been introduced in \cite{communi} and has been used to give alternative proofs of interesting results, such as the Reduction Theorem and the simplicity criteria for Leavitt path algebras (see \cite{communi, gabriela, reduction}) and is related to branching systems, see \cite{canto}.

 Ultragraphs are versatile objects that generalize graphs. In the algebraic setting, their associated algebras include Leavitt path algebras and algebras associated with infinite matrices and, in the analytical setting, their associated algebras generalized graph C*-algebras and Exel-Laca algebras, see \cite{Tom}. As with Leavitt path algebras of graphs and graph C*-algebras, the algebras associated with ultragraphs have deep connections with dynamical systems, see \cite{GGW, Tasca} for example. Furthermore, ultragraphs have applications in computer science, see \cite{computer}.

Our goal in this paper is to extend the results regarding the strong $\Z$-grading and epsilon-strong $\Z$ grading mentioned above to ultragraph Leavitt path algebras (which also have a $\Z$-grading, see \cite{Imanfar}), as well as find characterizations of ultragraph Leavitt path algebras that are strongly graded and epsilon-strong graded with respect to the free group on the edged grading (such grading was introduced in \cite{australia}). Our results establish a connection between the combinatorial structure of an ultragraph and its algebraic structure. Furthermore, the results regarding the grading over the free group of the edges are new even in the context of graphs.

We now give a detailed overview of the paper. In the next section, we present some necessary preliminaries and extend Condition~(Y) to ultragraphs. In Section~3, we show that the ultragraph Leavitt path algebra associated with an ultragraph $\GG$, $L_R(\GG)$, is strongly $\Z$-graded if, and only if, no vertex of $\mathcal{G}$ is a sink, $\mathcal{G}$ is row-finite, and $\mathcal{G}$ satisfies Condition~(Y), see Theorem~\ref{siri}. Aiming at a parallel with the analytical theory, we show in Section~4 that $L_R(\GG)$ is strongly graded if, and only if, the gauge action on the ultragraph C*-algebra, $C^*(\GG)$, is saturated. Given an ultragraph $\GG$, in Section~5 we give a sufficient condition for $L_R(\GG)$ to be epsilon-strongly $\Z$-graded  and show that this condition is not necessary. We also provide a necessary condition for $L_R(\GG)$ to be epsilon-strongly $\Z$-graded, see Proposition~\ref{tutucao}. Denote the free group of the edges on an ultragraph $\GG$ by $\F$. In Section~6, we recall the $\F$ grading on $L_R(\GG)$, show that $L_R(\GG)$ is strongly-$\F$-graded if, and only if, $\GG$ has only one edge $e$, in which case the source and range of $e$ are equal (see Theorem~\ref{voandonaagua}), and show that $L_R(\GG)$ is epsilon-strongly-$\F$-graded if, and onlly if, $L_R(\GG)$ is unital, what happens if, and only if, the set of vertices of $\GG$ is a generalized vertex of $\GG$ (see Theorem~\ref{ornitorrinco} and Corollary~\ref{platapus}).

\section{Preliminaries and Condition (Y)}

In this section we briefly recall the necessary notions used in this paper, set up notation, and introduce Condition~(Y) for ultragraphs. We start with the definition of an ultragraph (as in \cite{Tom}).

\begin{definicao}\label{def of ultragraph}
An \emph{ultragraph} is a quadruple $\mathcal{G}=(G^0, \mathcal{G}^1, r,s)$ consisting of two countable sets $G^0, \mathcal{G}^1$, a map $s:\mathcal{G}^1 \to G^0$, and a map $r:\mathcal{G}^1 \to P(G^0)\setminus \{\emptyset\}$, where $P(G^0)$ stands for the power set of $G^0$.
\end{definicao}

\begin{definicao}\label{def of mathcal{G}^0}
Let $\mathcal{G}$ be an ultragraph. Define $\mathcal{G}^0$ to be the smallest subset of $P(G^0)$ that contains $\{v\}$ for all $v\in G^0$, contains $r(e)$ for all $e\in \mathcal{G}^1$, and is closed under finite unions and non-empty finite intersections. Elements of $\mathcal{G}^0$ are called generalized vertices.
\end{definicao}


Next, we recall the definition of the Leavitt path algebra associated with an ultragraph. We remark that there are two similar ways to define the algebra (see \cite{australia, Imanfar}), depending on the definition of the generalized vertices, which have been shown in \cite{gildada} to coincide. We use the definition given in \cite{australia}.

\begin{definicao}\label{def of ultragraph algebra}
Let $\mathcal{G}$ be an ultragraph and $R$ be a unital commutative ring. The Leavitt path algebra of $\mathcal{G}$, denoted by $L_R(\mathcal{G})$, is the universal algebra with generators $\{s_e,s_e^*:e\in \mathcal{G}^1\}\cup\{p_A:A\in \mathcal{G}^0\}$ and relations
\begin{enumerate}
\item $p_\emptyset=0,  p_Ap_B=p_{A\cap B},  p_{A\cup B}=p_A+p_B-p_{A\cap B}$, for all $A,B\in \mathcal{G}^0$;
\item $p_{s(e)}s_e=s_ep_{r(e)}=s_e$ and $p_{r(e)}s_e^*=s_e^*p_{s(e)}=s_e^*$ for each $e\in \mathcal{G}^1$;
\item $s_e^*s_f=\delta_{e,f}p_{r(e)}$ for all $e,f\in \mathcal{G}$;
\item $p_v=\sum\limits_{s(e)=v}s_es_e^*$ whenever $0<\vert s^{-1}(v)\vert< \infty$.
\end{enumerate}
\end{definicao}
Let $\mathcal{G}$ be an ultragraph. A finite path is either an element of $\mathcal{G}^0$ or a sequence of edges $e_1...e_n$, with length $|e_1...e_n|=n$, and such that $s(e_{i+1})\in r(e_i)$ for each $i\in \{0,...,n-1\}$. An infinite path is a sequence $e_1e_2e_3...$, with length $|e_1e_2...|=\infty$, such that $s(e_{i+1})\in r(e_i)$ for each $i\geq 0$. The set of finite paths in $\mathcal{G}$ is denoted by $\mathcal{G}^*$, and the set of infinite paths in $\mathcal{G}$ is denoted by $\mathfrak{p}^\infty$. We extend the source and range maps as follows: $r(\alpha)=r(\alpha_{|\alpha|})$, $s(\alpha)=s(\alpha_1)$ for $\alpha\in \mathcal{G}^*$ with $0<|\alpha|<\infty$, $s(\alpha)=s(\alpha_1)$ for each $\alpha\in \mathfrak{p}^\infty$, and $r(A)=A=s(A)$ for each $A\in \mathcal{G}^0$. An element $v\in G^0$ is a sink if $s^{-1}(v) = \emptyset$, and we denote the set of sinks in $G^0$ by $G^0_s$. We say that $A\in \mathcal{G}^0$ is a sink if each vertex in $A$ is a sink. A vertex $v\in G^0$ is a source if $v\notin r(e)$ for any edge $e$, and $v$ is an infinite emitter if $|s^{-1}(v)|=\infty$. Finally, we denote the shift map in $\mathfrak{p}^\infty$ by $\sigma$, that is, $\sigma:\mathfrak{p}^\infty \rightarrow \mathfrak{p}^\infty$ is given by $\sigma(e_1e_2e_3\ldots)=e_2e_3\ldots$.

\subsection{Condition~(Y) }

In this subsection, we introduce Condition~(Y) for ultragraphs. This generalizes the definition of Condition~(Y) given for graphs in \cite{Clark} and recently studied in \cite{strongradedleavitt}.

\begin{definicao}\label{conditionY}
An ultragraph $\mathcal{G}$ satisfies Condition~(Y) if for each infinite path $e_1e_2e_3...\in \mathfrak{p}^\infty$ there exists a finite path $\alpha\in \mathcal{G}^*$ and some $k\in \N$ such that $s(e_{k+1})\in r(\alpha)$ and $|\alpha|=k+1$, that is, such that $\alpha \sigma^k(e_1e_2e_3\ldots) \in \mathfrak{p}^\infty$. 
\end{definicao}

Informally speaking, an ultragraph satisfies Condition~(Y) if, for each infinite path $p$ of this ultragraph, it is possible to replace some initial sub-path of length $k$ of $p$ by another finite path with length $k+1$, and still obtain an infinite path. 

\begin{exemplo}\label{ex2}
Let $\mathcal{G}$ be the following ultragraph:

\vspace{2cm}

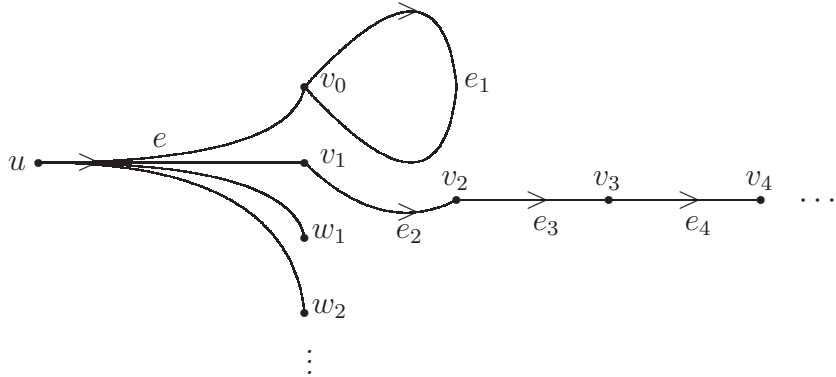
\centerline{
\setlength{\unitlength}{1cm}
\begin{picture}(9,0)
\put(-0.5,0){\circle*{0.1}}
\put(-0.9,-0.1){$u$}
\qbezier(-0.5,0)(2.8,0)(3,1)
\qbezier(-0.5,0)(2.8,0)(3,0)
\qbezier(-0.5,0)(2.8,0)(3,-1)
\qbezier(-0.5,0)(2.8,0)(3,-2)
\put(3,1){\circle*{0.1}}
\put(3,0){\circle*{0.1}}
\put(3,-1){\circle*{0.1}}
\put(3,-2){\circle*{0.1}}
\put(3,-2.8){\vdots}
\qbezier(3,1)(4.8,3)(5,1)
\qbezier(3,1)(4.8,-1)(5,1)
\qbezier(3,0)(4,-1)(5,-0.5)
\put(5,-0.5){\circle*{0.1}}
\qbezier(5,-0.5)(6,-0.5)(9,-0.5)
\put(7,-0.5){\circle*{0.1}}
\put(9,-0.5){\circle*{0.1}}
\put(4.8,-0.3){$v_2$}
\put(6.8,-0.3){$v_3$}
\put(8.8,-0.3){$v_4$}
\put(3.2,1){$v_0$}
\put(3.2,0){$v_1$}
\put(3.1,-1){$w_1$}
\put(3.1,-2){$w_2$}
\put(1,0.2){$e$}
\put(5.1,1){$e_1$}
\put(0,-0.1){$>$}
\put(4.2,1.9){$>$}
\put(4.2,-0.76){$>$}
\put(4.2,-1){$e_2$}
\put(5.9,-0.6){$>$}
\put(6,-0.9){$e_3$}
\put(7.9,-0.6){$>$}
\put(8,-0.9){$e_4$}
\put(9.5,-0.6){$\cdots$}
\end{picture}}
\vspace{3 cm}
\end{exemplo}

Notice that the ultragraph above does not satisfy Condition~(Y), since for the infinite path $p=ee_2e_3e_4...$ it is not possible to obtain an infinite path by replacing some initial sub-path $\beta$ of $p$ by another finite path $\alpha$ with $|\alpha|=|\beta|+1$.

\begin{remark}\label{condYoriginal} 
It is straightforward to check that an ultragraph satisfies Condition~(Y) if, and only if, for each infinite path $p\in \mathfrak{p}^\infty$ and for each $m\in \N$, it is possible to replace a initial sub-path $\beta$ of $p$ by another finite path $\alpha$, with $|\alpha|=|\beta|+m$, and still get an infinite path. In more symbolic notation, an ultragraph satisfies Condition~(Y) if, and only if, for every $p\in \mathfrak{p}^\infty$ and $m\in \N$, there exists $\alpha \in \mathcal{G}^{*}$, with $|\alpha|=|\beta|+m$, such that $\alpha \sigma^{|\beta|}(p)\in \mathfrak{p}^\infty$.
\end{remark}

\begin{remark}
In the graph setting, the original definition of Condition~(Y) (see \cite{Clark}) was the one described in Remark~\ref{condYoriginal}. Definition~\ref{conditionY} was introduced in \cite{strongradedleavitt}, where it was called Condition~(Y1) and where it was proved to be equivalent to the Condition~(Y) of \cite{Clark}. Since both descriptions are also equivalent in the ultragraph setting, we will not distinguish between Condition~(Y) and (Y1) in this paper and will refer to it only by Condition~(Y).
\end{remark}

 We finish this short subsection observing that if $\mathcal{G}$ is an ultragraph without sources then $\mathcal{G}$ satisfies Condition~(Y), but the converse is not true, as we illustrate in the next example.
 \begin{exemplo}
 Let $\mathcal{G}$ be the graph with two edges $e$ and $f$, and two vertices, $u$ and $v$, such that $s(e)=u$, $r(e)=v$, $s(f)=v=r(f)$. In this case, $u$ is a source of $\mathcal{G}$, but $\mathcal{G}$ satisfies Condition~(Y).  
\end{exemplo}

\section{The strong $\Z$-grading in $L_R(\mathcal{G})$}

In this section, we will describe the ultragraphs for which the associated Leavitt path algebra is strongly $\Z-$graded. We start recalling the definition of a graded ring and the definition of the $\Z-$grading on an ultragraph Leavitt path algebra.

\begin{definicao}
Let $T$ be a ring and $F$ be a group. We say that $T$ is $F$-graded if there is a collection of additive subgroups $\{ T_n \}_{n \in F}$ of $T$ with the following two properties.
\begin{enumerate}
\item $T = \bigoplus_{n \in F}T_n$.
\item $T_m T_n \subseteq T_{mn}$ for all $m,n \in F$.
\end{enumerate}
In this case, we also say that $\{ T_n \}_{n \in F}$ is a grading on $T$.
We call the subgroup $T_n$ the homogeneous component of $T$ of degree $n$, and we say that the elements of $T_n$ are homogeneous elements of degree $n$. 
If, in a graded ring $T$, the condition $T_mT_n=T_{mn}$ is satisfied for each $m,n\in F$, we say that $T$ is strongly $F$-graded.
\end{definicao}

If $\mathcal{G}$ is an ultragraph, then we may define a $\Z$-grading (\cite[Theorem~2.9]{Imanfar}) on the associated Leavitt path algebra $L_R(\mathcal{G})$ by setting, for each $m\in \Z$, 
\[L_R(\mathcal{G})_m=\left \{\sum\limits_{i=1}^k \lambda_is_{\alpha_i}p_{A_i}s_{\beta_i}^*:k\in \N, \alpha_i, \beta_i\in \mathcal{G}^*, A_i\in \mathcal{G}^0, \lambda_i\in R \text{ and } |\alpha_i|-|\beta_i|=m\right \}.\]



For future use, we record now an auxiliary result.

\begin{lema}\label{lemaauxiliar}

Let $\GG$ be an ultragraph and $\{L_R(\mathcal{G})_n\}_{n\in \Z}$ be the $\Z$-grading on $L_R(\mathcal{G})$ as above. For each $n\in \Z$, the following equalities hold: $$L_R(\mathcal{G})_nL_R(\mathcal{G})_0=L_R(\mathcal{G})_n=L_R(\mathcal{G})_0L_R(\mathcal{G})_n.$$
\end{lema}

\demo We show that  $L_R(\mathcal{G})_nL_R(\mathcal{G})_0=L_R(\mathcal{G})_n$ for each $n\in \Z$, and leave the proof that $L_R(\mathcal{G})_n=L_R(\mathcal{G})_0L_R(\mathcal{G})_n$ to the reader, since it is analogous. 

Notice that we only need to show that $L_R(\mathcal{G})_n\subseteq L_R(\mathcal{G})_nL_R(\mathcal{G})_0$ since, by the grading, $L_R(\mathcal{G})_nL_R(\mathcal{G})_0\subseteq L_R(\mathcal{G})_n$.

Let $x= \sum\limits_{i=1}^k\lambda_is_{\alpha_i}p_{A_i}s_{\beta_i}^*\in L_R(\mathcal{G})_n$. Then, for each $i$, there exists $B_i\in \mathcal{G}^0$ such that $s_{\alpha_i}p_{A_i}s_{\beta_i}^*p_{B_i}=s_{\alpha_i}p_{A_i}s_{\beta_i}^*$. Let $B=\bigcup\limits_{i=1}^kB_i$. Then, $p_B\in L_R(\mathcal{G})_0$ and $p_{B_i}p_B=p_{B_i}$. Hence  $$x=\sum\limits_{i=1}^k\lambda_is_{\alpha_i}p_{A_i}s_{\beta_i}^*=\sum\limits_{i=1}^k\lambda_is_{\alpha_i}p_{A_i}s_{\beta_i}^*p_{B_i}=(\sum\limits_{i=1}^k\lambda_is_{\alpha_i}p_{A_i}s_{\beta_i}^*p_{B_i})p_B\in L_R(\mathcal{G})_nL_R(\mathcal{G})_0,$$
and therefore $L_R(\mathcal{G})_n\subseteq L_R(\mathcal{G})_nL_R(\mathcal{G})_0$, as desired.  
\fim

It is shown in \cite{strongradedleavitt} that if $E$ is a directed graph then $L_R(E)$, with its canonical $\Z$-grading, is strongly $\Z$-graded  if, and only if, $E$ has no sinks, is row-finite, and satisfies Condition~(Y). We will show that the same holds for ultragraphs, but for this, we need a version of row-finitess in the ultragraph context. Searching the literature, we have found the following concept, which fits nicely in our work.

\begin{definicao}\cite[Definition~2.1]{BP1} An ultragraph $\mathcal{G}$ is row-finite if $r(e)$ is finite for every $e\in \mathcal G^1$ and there are no infinite infinite emitters in $G^0$.
\end{definicao}

\begin{exemplo} Let $\mathcal{G}$ be the ultragraph as follows: 

\vspace{1cm}
\centerline{
\setlength{\unitlength}{1cm}
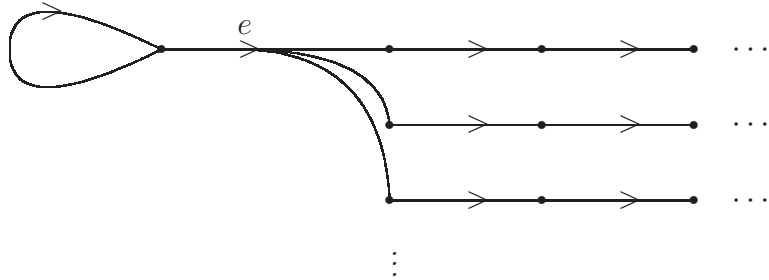
\begin{picture}(5,0)
\put(0,0){\circle*{0.1}}
\qbezier(0,0)(2,0)(7,0)
\qbezier(1,0)(2.9,0)(3,-1)
\qbezier(1,0)(2.9,0)(3,-2)
\put(3,0){\circle*{0.1}}
\put(3,-1){\circle*{0.1}}
\put(3,-2){\circle*{0.1}}
\qbezier(0,0)(-2,-1)(-2,0)
\qbezier(0,0)(-2,1)(-2,0)
\put(-1.6,0.4){$>$}
\put(1,-0.1){$>$}
\put(1,0.2){$e$}
\qbezier(3,-1)(4,-1)(7,-1)
\qbezier(3,-2)(4,-2)(7,-2)
\put(5,-1){\circle*{0.1}}
\put(7,-1){\circle*{0.1}}
\put(5,-2){\circle*{0.1}}
\put(7,-2){\circle*{0.1}}
\put(5,0){\circle*{0.1}}
\put(7,0){\circle*{0.1}}
\put(7.5,-0.1){$\cdots$}
\put(7.5,-1.1){$\cdots$}
\put(7.5,-2.1){$\cdots$}
\put(3,-3){$\vdots$}
\put(4,-0.1){$>$}
\put(6,-0.1){$>$}
\put(4,-1.1){$>$}
\put(6,-1.1){$>$}
\put(4,-2.1){$>$}
\put(6,-2.1){$>$}
\end{picture}}
\vspace{4 cm}

This ultragraph is not row-finite, since $r(e)$ is an infinite set. It has no sinks (and no sources) and satisfies Condition~(Y). Due to the following lemma, $L_R(\mathcal{G})$ is not strongly $\Z$-graded. 

\end{exemplo}

\begin{lema}\label{r(e)finite}
Let $\mathcal{G}$ be an ultragraph such that $r(e)$ is infinite for some edge $e$. Then $L_R(\mathcal{G})$ is not strongly $\Z$-graded.
\end{lema}

\demo Let $e\in \mathcal{G}^1$ be such that $r(e)$ is infinite. Suppose that $L_R(\mathcal{G})$ is strongly $\Z$-graded. Then, $p_{r(e)}\in L_R(\mathcal{G})_1L_R(\mathcal{G})_{-1}$, and hence $p_{r(e)}=\sum\limits_{k=1}^m a_kb_k$ with $a_k\in L_R(\mathcal{G})_1$ and $b_k\in L_R(\mathcal{G})_{-1}$. Now, for each $k$, write $$a_k=\sum\limits_{i\in F_k}\lambda_is_{\alpha_i}p_{A_i}s_{\beta_i}^*,$$ where $F_k\subseteq \N$ is a finite set, $\lambda_i\in R$,  $A_i\in \mathcal{G}^0$, and $\alpha_i$, $\beta_i$ are finite paths in $\mathcal{G}$ with $|\alpha_i|\geq 1$ for each $i$.  Since $r(e)$ is infinite, there exists $v\in r(e)$ such that $v\neq s(\alpha_i)$ for each $i\in F_k$ and each $k\in \{1,...,m\}$. Then $p_va_k=0$ for each $k$, and so $$p_v=p_vp_{r(e)}=p_v\sum\limits_{k=1}^m a_kb_k=0$$ which is impossible, since $p_v\neq 0$. So, it follows that $L_R(\mathcal{G})$ is not strongly $\Z$-graded.  
\fim

The above lemma tells us that strongly graded ultragraph Leavitt path algebras only occur associated with ultragraphs for which the range of each edge is finite. We will use this class of ultragraphs a sufficient number of times to give it a name, see below.

\begin{definicao} We say that an ultragraph is \emph{finite range} if each edge has finite range.
\end{definicao}

We show next that each algebra in the class of Leavitt path algebras associated with finite range ultragraphs can be realized as the Leavitt path algebra of a graph (This is undoubtedly known to experts, we include here for completeness). To do this, we first describe how to associate a graph with an ultragraph.



Let $\mathcal{G}$ be an ultragraph, and let $E_{\mathcal G}=(E_{\mathcal G}^0, E_{\mathcal G}^1, r,s)$ be the directed graph associated to $\mathcal{G}$ as follows: define $E_{\mathcal G}^0=G^0$, $E_{\mathcal G}^1=\{e_u:e\in \mathcal{G}^1 \text{ and }u\in r(e)\}$, and define the range and source maps in $E_{\mathcal G}$ by $r(e_u)=u$ and $s(e_u)=s_1(e)$. 

\begin{proposicao}\label{isographandultragraph}
Let $\mathcal{G}$ be an ultragraph such that $r(e)$ is finite for each edge $e$ and let $E_{\mathcal G}$ be the directed graph associated to $\mathcal{G}$ as above. Then there exists a $\Z$-graded isomorphism $\pi:L_R(\mathcal{G})\rightarrow L_R(E_{\mathcal G})$, with respect to the natural $\Z$-gradings of $L_R(\mathcal{G})$ and $L_R(E_{\mathcal G})$.
\end{proposicao}

\demo Since, for each edge $e$ of $\mathcal G$, $r(e)$ is a finite set, we have that each $A\in \mathcal{G}^0$ is also finite. From the universal property of $L_R(\mathcal{G})$, we get a homomorphism $\pi:L_R(\mathcal{G})\rightarrow L_R(E_{\mathcal G})$ such that $\pi(p_A)=\sum\limits_{v\in A}p_v$, for each $A\in \mathcal{G}^0$, and $\pi(s_e)=\sum\limits_{u\in r(e)}s_{e_u}$, for each $e\in \mathcal{G}^1$. Similarly, there exists a homomorphism $\varphi:L_R(E_{\mathcal G})\rightarrow L_R(\mathcal{G})$ such that $\varphi(s_{e_u})=s_ep_u$, for each $e_u\in E_{\mathcal G}^1$, and $\varphi(p_u)=p_u$, for each $u\in E_{\mathcal G}^0$. Notice that $\pi$ and  $\varphi$ are inverses of each other. Now, since $\pi$ is a homomorphism, it is easy to see that $\pi(L_R(\mathcal{G})_m)\subseteq L_R(E_{\mathcal G})_m$, and also that $\pi^{-1}(L_R(E_{\mathcal G})_m)\subseteq L_R(\mathcal{G})_m$ for each $m\in \Z$. Therefore, $\pi(L_R(\mathcal{G})_m)= L_R(E_{\mathcal G})_m$ for each $m\in \Z$.
\fim

For general ultragraphs, we relate Condition~(Y) in $\mathcal G$ with Condition~(Y) in $E_{\mathcal G}$ in the result below.

\begin{proposicao}\label{Y-condition}
Let $\mathcal{G}$ be an ultragraph and $E_{\mathcal G}$ be the associated directed graph as above. Then $\mathcal{G}$ satisfies Condition~(Y) if, and only if, $E_{\mathcal G}$ satisfies Condition~(Y).
\end{proposicao}

\demo Suppose that $\mathcal{G}$ satisfies Condition~(Y). Let ${e_1}_{u_1}{e_2}_{u_2}{e_3}_{u_3}...$ be an infinite path in $E_{\mathcal G}$, where each ${e_i}_{u_i}\in E_{\mathcal G}^1$. Then $e_1e_2e_3....$ is an infinite path of $\mathcal{G}$. Since $\mathcal{G}$ satisfies Condition~(Y), there exists $k\in \N$ and a finite path $x_1....x_{k+1}$ in $\mathcal{G}$, with $x_i\in \mathcal{G}^1$ for each $i$, such that $x_1x_2...x_{k+1}e_{k+1}e_{k+2}...$ is an infinite path of $\mathcal{G}$. Let $v_i$ be the source of $x_{i+1}$ for each $i\in \{1,...,k\}$, and let $v_{k+1}$ be the source of $e_{k+1}$ (in $\mathcal{G}$). Then  $${x_1}_{v_1}{x_2}_{v_2}...{x_{k+1}}_{v_{k+1}}{e_{k+1}}_{u_{k+1}}{e_{k+2}}_{u_{k+2}}...$$ is an infinite path of $E_{\mathcal G}$. Therefore $E_{\mathcal G}$ satisfies Condition~(Y).

Similarly one shows that if $E_{\mathcal G}$ satisfies Condition~(Y) then $\mathcal{G}$ also satisfies Condition~(Y).
\fim

The main theorem of this section is the following characterization of strongly graded ultragraph Leavitt path algebras.

\begin{teorema}\label{siri}
Let $\mathcal{G}$ be an ultragraph and $R$ a unital ring. Then $L_R(\mathcal{G})$ is strongly $\Z$-graded, with the grading $\{L_R(\mathcal{G})_n\}_{n\in \Z }$, if and only if no vertex of $\mathcal{G}$ is a sink, $\mathcal{G}$ is row-finite, and $\mathcal{G}$ satisfies Condition~(Y). 
\end{teorema}

\demo
Suppose that $L_R(\mathcal{G})$ is strongly $\Z$-graded. By Lemma~\ref{r(e)finite}, $\mathcal{G}$ is a finite range ultragraph. Let $E_{\mathcal G}$ be the directed graph associated to $\mathcal{G}$. By Proposition~\ref{isographandultragraph}, $L_R(E_{\mathcal G})$ is also strongly $\Z$-graded. It follows from \cite[Lemma 3.1]{strongradedleavitt} that $E_{\mathcal G}$ has no sinks, is row-finite and satisfies Condition~(Y). This implies that $\mathcal{G}$ is also row-finite and has no sinks. Moreover, since $E_{\mathcal G}$ satisfies Condition~(Y), by Proposition~\ref{Y-condition},  we conclude that $\mathcal{G}$ also satisfies Condition~(Y).

Now, suppose that $\mathcal{G}$ satisfies Condition~(Y), is row-finite and no vertex of $\mathcal{G}$ is a sink. Then, the same holds for $E_{\mathcal G}$, that is, $E_{\mathcal G}$ has no sinks, is row-finite, and (by Proposition~\ref{Y-condition}) satisfies Condition~(Y). Therefore, \cite[Proposition 4.3]{strongradedleavitt} implies that  $L_R(E_{\mathcal G})$ is strongly $\Z$-graded. Since $\mathcal{G}$ is finite range, by Proposition~\ref{isographandultragraph}, we conclude that $L_R(\mathcal{G})$ is also strongly $\Z$-graded.
\fim

\begin{remark}
The above theorem can also be obtained using the realization of ultragraph Leavitt path algebras as Steinberg algebras given in \cite{gildada} and applying the results of \cite{Clark} (in particular Theorem~3.11 in \cite{Clark}). In fact, for ultragraphs withouth sinks, the above result was proved via groupoids in \cite[Corollary~3.4]{HaNam}. However, there is an imprecison in the groupoid description of an ultragraph Leavitt path algebra given in \cite{HaNam}, which was explaind in \cite[Remark~3.7]{gildada}. All considered, we believe that our direct approach is simpler, as it avoids the intricate topology used in \cite{gildada}.
\end{remark}

\section{Gauge action on C*($\mathcal G$) and Condition~(Y)}

It is well known that Leavitt path algebras have strong links with its analytical counter-parts, graph C*-algebras. What is still mysterious, in many settings, is to determine when a property of the Leavitt path algebra associated to a graph implies a property of the graph C*-algebra associated to this graph, and vice-versa (see \cite[Chapter~5]{book}). With ultragraphs algebras, the same kind of phenomena occurs. In this section, we show that, for ultragraphs, a strong $\Z$-grading on an ultragraph Leavitt path algebra is equivalent to the saturation of the gauge action on the corresponding ultragraph C*-algebra. This extends to ultragraphs the results for graphs proved in \cite{chirvasitu}. Furthermore, using what was previously proved for ultragraph Leavitt path algebras (Theorem~\ref{siri}), we get a characterization of ultragraph C*-algebras with saturated gauge actions. We start recalling some of the key concepts and setting up notation.

Let $\mathcal{G}$ be an ultragraph,  $C^*(\mathcal{G})$ be the associated ultragraph $C^*$-algebra (for details, see \cite{Tom}), and $S^1=\{z\in \C:|z|=1\}$. The \emph{gauge action} associated to $\mathcal G$ is the map from $S^1$ to the automorphisms of $C^*(\mathcal G)$ given by $S^1\ni z\mapsto \gamma_z\in Aut(C^*(\mathcal{G}))$, where $\gamma_z(s_e)=zs_e$ for each $e\in \mathcal{G}^1$, and $\gamma_z(p_A)=p_A$ for each $A\in \mathcal{G}^0$. Next we define saturated gauge actions. Notice that in \cite{chirvasitu} the author refers to saturated gauge actions as free gauge actions. We will follow the nomenclature defined in \cite{chris}.


\begin{definicao}\cite[Definition 1.8]{chirvasitu} Let $G$ be a compact abelian group acting on a $C^*$-algebra $A$. For every character $\tau$ of $G$ (that is, $\tau$ is a homomorphism from $G$ to $S^1$), denote by $A_\tau$ the $\tau$-eigenspace of the action, that is,
$$A_\tau=\{a\in A:ga=\tau(g)a\,\,\,\,\forall{g\in G}\}.$$
 The fixed point algebra is the eigenspace associated to the constant character 1 and is denoted by $A^G$. Finally, the gauge action is saturated if $A_\tau^*A_\tau$ is dense in the fixed point algebra $A_1=A^G$, for each character $\tau$ of $G$.
\end{definicao}

\begin{remark}
Recall that the character space $\widehat{S^1}$ of $S^1$ is $\Z$. So, given $\tau \in \widehat{S^1}$, there exists $n\in \Z$ such that $\tau(z)=z^{n}$ for each $z\in S^1$. Therefore, the $\tau$-eigenspace of the gauge action $\gamma$ of an ultragraph $C^*$-algebra $C^*(\mathcal{G})$ is the space $$C^*(\mathcal{G})_\tau=\{a\in C^*(\mathcal{G}): \gamma_z(a)=z^na\,\,\,\, \forall{z\in S^1}\}.$$ From now on, we denote the space $C^*(\mathcal{G})_\tau$ by $C^*(\mathcal{G})_n$. With this notation, the fixed point algebra is denoted by  $C^*(\mathcal{G})_0$.
\end{remark}

Let $\GG$ be an ultragraph, $L_\C(\mathcal{G})$ be the associated ultragraph Leavitt path algebra over the complex numbers, and let $\{L_\C(\mathcal{G})_m\}_{m\in \Z}$ be the associated $\Z$-grading. It follows directly from the definition of $L_\C(\mathcal{G})_m$ that $L_\C(\mathcal{G})_m\subseteq C^*(\mathcal{G})_m$ for each $m\in \Z$. Moreover, using the conditional expectation associated with the gauge action, it follows that $L_\C(\mathcal{G})_0$ is dense in the fixed point algebra $C^*(\mathcal{G})_0$.

We now have the necessary elements to prove the main result of this section. In fact, since we have already built the necessary auxiliary results, the outline of the proof follows analogously to the proof of \cite[Proposition~2.13]{chirvasitu}. We include it here for completeness.

\begin{teorema}
Let $\mathcal{G}$ be an ultragraph. Then the gauge action $\gamma$ on $C^*(\mathcal{G})$ is saturated if, an only if, $L_\C(\mathcal{G})$ is strongly $\Z$-graded, with respect to the grading $\{L_\C(\mathcal{G})_m\}_{m\in \Z}$.
\end{teorema}

\demo Suppose that $L_\C(\mathcal{G})$ is strongly $\Z$-graded. 
Then, for each $m\in \Z$, we have that $C^*(\mathcal{G})_{-m}C^*(\mathcal{G})_m \supseteq L_\C(\mathcal{G})_{-m}L_\C(\mathcal{G})_m = L_\C(\mathcal{G})_0$. Since $L_\C(\mathcal{G})_0$ is dense in $C^*(\mathcal{G})_0$, it follows that $C^*(\mathcal{G})_{-m}C^*(\mathcal{G})_m$ is also dense  $C^*(\mathcal{G})_0$. 


For the converse, suppose that the gauge action $\gamma$ is saturated. To show that $L_\C(\mathcal{G})$ is strongly $\Z$-graded, by \cite[Proposition 2.1]{strongradedleavitt} and Lemma~\ref{lemaauxiliar}, it is enough to prove that $L_\C(\mathcal{G})_1L_\C(\mathcal{G})_{-1}=L_\C(\mathcal{G})_0=L_\C(\mathcal{G})_{-1}L_\C(\mathcal{G})_1$.

We show that $L_\C(\mathcal{G})_1L_\C(\mathcal{G})_{-1}=L_\C(\mathcal{G})_0$ and leave the other case to the reader, as it is analogous. It is clear that $L_\C(\mathcal{G})_1L_\C(\mathcal{G})_{-1}\subseteq L_\C(\mathcal{G})_0$ and so we only have to prove  that $L_\C(\mathcal{G})_0\subseteq L_\C(\mathcal{G})_1L_\C(\mathcal{G})_{-1}$. For this, it is enough to show that $s_\mu p_A s_\nu^*\in L_\C(\mathcal{G})_1L_\C(\mathcal{G})_{-1}$, for each $A\in \mathcal{G}^0$ and finite paths $\mu,\nu$ with $|\mu|=|\nu|$.
Furthermore, since $s_\mu s_\mu^*$ is a unit on the left for $s_\mu p_A s_\nu^*$ and $L_\C(\mathcal{G})_1L_\C(\mathcal{G})_{-1}L_\C(\mathcal{G})_0\subseteq L_\C(\mathcal{G})_1L_\C(\mathcal{G})_{-1}$, it is sufficient to prove that $s_\mu s_\mu^*\in L_\C(\mathcal{G})_1L_\C(\mathcal{G})_{-1}$. Let $x=s_\mu s_\mu^*$.
Since $\gamma$ is saturated, there are finitely many $x_i,y_i\in L_\C(\mathcal{G})_1$ such that $$||\sum\limits_ix_iy_i^*-x||<1.$$ 
 Let $a_i=xx_i$ and $b_i=xy_i$. Since $||x||=1$ and $x^2=x=x^*$, we get that $$||\sum\limits_i a_i b_i^*-x||=||\sum\limits_ixx_iy_i^*x^*-x||=||x(\sum\limits_ix_iy_i^*-x)x^*||<1.$$ Notice that $x$ is a unit for all the elements $a_ib_i^*$. Let $A$ be the finite dimensional C*-algebra generated by all $a_i b_i^*$ and $x$, and let $B$ be the finite dimensional C*-algebra generated by all $a_i b_i^*$. Since $$||\sum\limits_i a_i b_i^*-x||<1$$ and $x$ is the unit of $A$, we have that $\sum\limits_i a_i b_i^*$ is invertible. Therefore, $B$ contains an invertible element and, since $B$ is an ideal of $A$, we obtain that $B=A$. Hence, $x\in B\subseteq L_\C(\mathcal{G})_1L_\C(\mathcal{G})_{-1}$ as desired.

\fim

Joining the above result with Theorem~\ref{siri}, we get the following characterization of ultragraph C*-algebras with saturated gauge action.

\begin{corolario} Let $\GG$ be an ultragraph. The following are equivalent.

\begin{itemize}
\item The gauge action on $C^*(\mathcal{G})$ is saturated.
\item $L_{\C}(\mathcal{G})$ is strongly $\Z$-graded. 
\item  $ \mathcal{G}$ has no sinks, is row-finite, and  satisfies Condition~(Y).
\end{itemize}

\end{corolario}

\begin{remark}
In the C*-algebraic setting, the above results generalizes \cite[Theorem~2.14]{chirvasitu} and \cite[Proposition~2]{won}. In the algebraic setting, the above result generalizes \cite[Theorem~4.2]{Clark}.
\end{remark}

\section{Epsilon strong $\Z$-grading on $L_R(\mathcal{G})$}

In an attempt to single out a class of well-behaved group
graded rings which naturally includes all crossed products by unital twisted partial actions,
Nystedt, Oinert, and Pinedo, see \cite{NOP},  recently introduced a generalization of unital strongly graded rings called epsilon-strongly graded rings. Among the properties with well described criteria for epsilon-strongly graded rings we find separability, semisimplicity, heredity, and Frobenius, see \cite{NOP}. In this section, we characterize ultragraph Leavitt path algebras that are epsilon-strongly $\Z$-graded in terms of the combinatorics of the ultragraph. 

Recall from \cite{groupgradleavitt} that a grading $\{T_g\}_{g\in G}$ of a ring $T$ is an {\it epsilon-strong $G$-grading} if, for each $g\in G$, the $(T_gT_{g^{-1}}\,-\,T_{g^{-1}}T_g)$ bimodule $T_g$ is unital. Moreover, $T$ is {\it symmetrically $G$-graded} if $T_gT_{g^{-1}}T_g=T_g$ for each $g\in G$.

The following result gives a useful description of epsilon-strongly $G-$graded rings.

\begin{proposicao}\label{epsilon-strongly}\cite[Proposition~7]{groupgradleavitt}
The following assertions are equivalent:
\begin{enumerate}
    \item The ring $T$ is epsilon-strongly $G$-graded, with epsilon-strong $G-$grading $\{T_g\}_{g\in G}$;
    \item The ring $T$ is symmetrically $G$-graded and, for each $g\in G$, the ring $T_gT_{g^{-1}}$ is unital;
    \item For each $g\in G$, there exists $\epsilon_g \in T_gT_{g^{-1}}$ such that, for all $s\in T_g$, the equalities $\epsilon_gs=s=s\epsilon_{g^{-1}}$ hold.
\end{enumerate}
\end{proposicao}

\begin{remark}\label{unitalremark} With the notation of the above proposition, notice that every epsilon-strongly $G$-graded ring $T$ is unital, with unit the element $\epsilon_e$, where $e$ is the unit of $G$. To verify this, let $g\in G$ and $s\in T_g$. Since $T$ is symmetrically $G$ graded, we can write $s=\sum\limits_{i=1}^ma_ib_ic_i,$ where $a_i,c_i\in T_g$ and $b_i\in T_{g^{-1}}$. Then $\epsilon_es=s$, since $a_ib_i\in T_e$, and $s\epsilon_e=s$, since  $b_ic_i\in T_e$.
\end{remark}

As we have seen above, the units play an important role in the study of epsilon-strongly $G$-graded rings. In the following lemma, we characterize unital ultragraph Leavitt path algebras. This is an algebraic version of the characterization of unital ultragraph C*-algebras given in \cite[Lemma~3.2]{Tom}.

\begin{lema}\label{unitallemma} (c.f. \cite[Lemma~3.2]{Tom}) Let $\mathcal G$ be an ultragraph. 
Te associated Leavitt path algebra $L_R(\mathcal{G})$ is unital if, and only if, $G^0\in \mathcal{G}^0$.
\end{lema}

\demo Suppose $G^0\in \mathcal{G}^0$. Notice that for finite paths $\alpha,\beta$, possibly with length zero, and for each element $A\in \mathcal{G}^0$ the equalities $p_{G^0}s_\alpha p_A s_\beta=s_\alpha p_A s_\beta=s_\alpha p_A s_\beta p_{G^0}$ are true. Since each element of $L_R(\mathcal{G})$ is a finite sum of elements of the form $s_\alpha p_A s_\beta$, it follows that $p_{G^0}$ is the unit of $L_R(\mathcal{G})$.

For the converse, suppose that $L_R(\mathcal{G})$ is unital, and denote its unit by $1$. Then $1= \sum\limits_{i=1}^m \lambda_i s_{\alpha_i} p_{A_i} s_{\beta_i}^*$, for some $\alpha_i,\beta_i \in \mathcal G^*, A_i \in \mathcal G^0$. 
Let $B= \bigcup\limits_{i=1}^m s(\alpha_i)\cup A_i \cup r(\beta_i)$. Then $B\in \mathcal G^0$. Suppose that $G^0 \notin \mathcal G^0$, and let $v \in G^0 \setminus B$. Then, \[p_v =  p_v 1 = p_v \sum\limits_{i=1}^m \lambda_i s_{\alpha_i} p_{A_i} s_{\beta_i}^* = 0,\] a contradiction.

\fim

As a consequence of the above results, we obtain that an ultragraph Leavitt path algebra that is epsilon-strongly $G$-graded always contains the set of vertices as a generalized vertex. We make this precise below.

\begin{proposicao}\label{G^0 unit}
Let $\mathcal{G}$ be an ultragraph and $R$ be a ring and $G$ be a group. If $L_R(\mathcal{G})$ is epsilon-strongly $G$-graded, then $G^0\in \mathcal{G}^0$.
\end{proposicao}

\demo The proof follows from Remark \ref{unitalremark} and Lemma \ref{unitallemma}.
\fim

It is proved in \cite[Theorem 28]{groupgradleavitt} that, for any finite directed graph $E$ and any group $G$, the Leavitt path algebra $L_R(E)$ is always epsilon-strongly $G$-graded. We give below a characterization of when the $\Z$-grading on an ultragraph Leavitt path algebra is epsilon-strong. This also provides a converse (when the group $G$ is $\Z$) for \cite[Theorem 28]{groupgradleavitt}, in the context of Leavitt path algebras, that is, our result implies that if a Leavitt path algebra associated to a graph is epsilon-strongly $\Z$-graded, then the graph is finite.

\begin{proposicao}\label{tutucao} Let $\mathcal{G}$ be an ultragraph. 
\begin{enumerate}
\item If $L_R(\mathcal{G})$ is epsilon-strongly $\Z$-graded, then $\mathcal{G}^1$ is a finite set and $G^0\in \mathcal{G}^0$.
\item If $\mathcal{G}^1$ is a finite set, $G^0\in \mathcal{G}^0$, and for each edge $e_1$ there exists an edge $e_2$ such that $s(e_1)\in r(e_2)$, then $L_R(\mathcal{G})$ is epsilon-strongly $\Z$-graded.
\end{enumerate}
\end{proposicao}

\demo Suppose that $L_R(\mathcal{G})$ is epsilon-strongly $\Z$-graded. It follows from Proposition~\ref{G^0 unit} that $G^0\in \mathcal{G}^0$. We show that $\mathcal{G}^1$ is finite. By Proposition~\ref{epsilon-strongly}, the ring  $L_R(\mathcal{G})_1L_R(\mathcal{G})_{-1}$ is unital. Let $x\in L_R(\mathcal{G})_1L_R(\mathcal{G})_{-1}$ be the unit. Then,  $x=\sum\limits_{i=1}^na_ib_i$ with $a_i\in L_R(\mathcal{G})_1$ and $b_i\in L_R(\mathcal{G})_{-1}$. Notice that each $a_i$ is of the form $a_i=\sum\limits_{j\in F_i}\lambda_js_{\alpha_j}p_{A_j}s_{\beta_j}^*$, where $\lambda_j\in R$, $\alpha_j,\beta_j\in \mathcal{G}^*$, $|\alpha_j|>0$, $A_j\in \mathcal{G}^0$ for all $j$, and $F_i$ is a finite set. Let $F$ be the set of all the first edges of all the $\alpha_j$, with $j\in F_i$ and $i\in \{1,...,n\}$, that is, \[F=\bigcup_{i=1}^n \{e_j:\alpha_j = e_j \alpha_j^{'}, j\in F_i\}.\] Clearly, $F$ is finite. Now, suppose that $\mathcal{G}^1$ is infinite. Chose $e\in \mathcal{G}^1\setminus F$. Then, $s_es_e^*\in L_R(\mathcal{G})_1L_R(\mathcal{G})_{-1}$ and $s_es_e^*x=0$ (since $s_e^*a_i=0$ for each $i$), which is impossible since $x$ is the unit of $L_R(\mathcal{G})_1L_R(\mathcal{G})_{-1}$. 
Therefore, $\mathcal{G}^1$ is a finite set.

Now, suppose that $\mathcal{G}^1$ is finite set, $G^0\in \mathcal{G}^0$, and that for each edge $e_1$ there exists an edge $e_2$ such that $s(e_1)\in r(e_2)$ . To show that $L_R(\mathcal{G})$ is epsilon-strongly $\Z$-graded, we will prove the second item of Proposition~\ref{epsilon-strongly}.

First notice that $L_R(\mathcal{G})_0 L_R(\mathcal{G})_0 =L_R(\mathcal{G})_0  $ is unital, with unit $p_{G^0}$.

Let $n\in \Z$ with $n\geq 1$. Let $F=\{\alpha \in \mathcal{G}^*:|\alpha|=n\}$, and notice that $F$ is a finite set. Let $x_n=\sum\limits_{\alpha\in F}s_{\alpha}s_{\alpha}^*$, and notice that $x_na=a$, for every $a\in L_R(\mathcal{G})_n$, and $bx_n=b$, for each $b\in L_R(\mathcal{G})_{-n}$. Therefore, $x_n$ is a unit of $L_R(\mathcal{G})_nL_R(\mathcal{G})_{-n}$ and, moreover, $L_R(\mathcal{G})_n\subseteq L_R(\mathcal{G})_nL_R(\mathcal{G})_{-n}L_R(\mathcal{G})_n$ (since we can write each $a\in L_R(\mathcal{G})_n$ as $a=x_na\in L_R(\mathcal{G})_nL_R(\mathcal{G})_{-n}L_R(\mathcal{G})_n$). Hence, 
\[L_R(\mathcal{G})_n= L_R(\mathcal{G})_nL_R(\mathcal{G})_{-n}L_R(\mathcal{G})_n.\]
Next, we show that $L_R(\mathcal{G})_{-n}L_R(\mathcal{G})_n$ is unital. Notice that, for each $\alpha\in F$, $s_{\alpha}^*s_{\alpha}\in L_R(\mathcal{G})_{-n} L_R(\mathcal{G})_n$. Moreover, if $\alpha, \beta\in F$, then $s_\alpha^*s_\alpha s_\beta^*s_\beta \in L_R(\mathcal{G})_{-n}L_R(\mathcal{G})_n$, since $s_\alpha^* \in L_R(\mathcal{G})_{-n}$ and $s_\alpha s_\beta^*s_\beta\in L_R(\mathcal{G})_nL_R(\mathcal{G})_0\subseteq L_R(\mathcal{G})_n$. 
From this we obtain that, for each $\alpha, \beta \in F$, the projection $p_{r(\alpha)\cup r(\beta)}$ belongs to $L_R(\mathcal{G})_{-n}L_R(\mathcal{G})_{n}$, since \[p_{r(\alpha)\cup r(\beta)}=p_{r(\alpha)}+p_{r(\beta)}-p_{r(\alpha)}p_{r(\beta)}= s_{\alpha}^*s_{\alpha}+s_{\beta}^*s_{\beta}-s_\alpha^*s_{\alpha}s_\beta^*s_{\beta}.\]

Let $A=\bigcup\limits_{\alpha\in F}r(\alpha)$. By induction, we obtain that $p_A \in L_R(\mathcal{G})_{-n}L_R(\mathcal{G})_n$ (since $F$ is finite.) Notice that, since $r(e)\subseteq A$ and $s(e)\in A$ for each edge $e$ then $p_A a=a$ for each $a\in L_R(\mathcal{G})_{-n}$ and $bp_A=b$ for each $b\in L_R(\mathcal{G})_n$, and so $p_A$ is a unit of $L_R(\mathcal{G})_{-n}L_R(\mathcal{G})_{n}$.

Finally, the proof that $L_R(\mathcal{G})_{-n}= L_R(\mathcal{G})_{-n}L_R(\mathcal{G})_{n}L_R(\mathcal{G})_{-n}$ is similar to what we did to prove that $L_R(\mathcal{G})_{n}= L_R(\mathcal{G})_{n}L_R(\mathcal{G})_{-n}L_R(\mathcal{G})_{n}$.

We have now checked that all conditions in the second item of Proposition~\ref{epsilon-strongly}, and hence it follows that $L_R(\mathcal{G})$ is epsilon-strongly $\Z$-graded.
\fim

\begin{remark} In the previous proposition, the hypothesis that for each edge $e_1$ there exists an edge $e_2$ such that $s(e_1)\in r(e_2)$ (besides $\mathcal{G}^1$ being finite and $G^0\in \mathcal{G}^0$) is sufficient to $L_R(\mathcal{G})$ be epsilon-strongly $\Z$-graded but is not necessary.   
For example, let $\mathcal{G}$ be the ultragraph with only one edge $e$, where $r(e)$ is an arbitrary set and $s(e)\notin r(e)$. Notice that $L_R(\mathcal{G})_n=0$ for each $|n|>1$ and that $L_R(\mathcal{G})_0L_R(\mathcal{G})_0=L_R(\mathcal{G})_0$ is unital, with unit $p_{r(e)\cup s(e)}$. Moreover, $L_R(\mathcal{G})_1L_R(\mathcal{G})_{-1}$ is unital with unit $p_{s(e)}=s_es_e^*$ and $L_R(\mathcal{G})_{-1}L_R(\mathcal{G})_{1}$ is unital with unit $p_{r(e)}=s_e^*s_e$. Finally it easy to see that $L_R(\mathcal{G})_1L_R(\mathcal{G})_{-1}L_R(\mathcal{G})_1=L_R(\mathcal{G})_1$  and  $L_R(\mathcal{G})_{-1}L_R(\mathcal{G})_1L_R(\mathcal{G})_{-1}=L_R(\mathcal{G})_{-1}.$ Now, by the second item of Proposition \ref{epsilon-strongly}, it follows that $L_R(\mathcal{G})$ is epsilon-strongly $\Z$-graded, but, for this ultragraph, it is not true that the source of each edge is contained in a range of some edge. 
\end{remark}

\section{The free group grading on $L_R(\mathcal{G})$}

In the previous sections, we studied in detail the $\Z$-grading on $L_R(\mathcal{G})$. Another interesting grading on $L_R(\mathcal{G})$ is determined by the free group of the edges, $\F$. This grading has been defined in \cite{australia}, generalizing previous work on Leavitt path algebras of graphs \cite{canto, communi, gabriela}, and it has been used to obtain interesting results on ultragraph (and graph) Leavitt path algebras, see \cite{canto, reduction, australia} for example. All the results in this section are also new for Leavitt path algebras of graphs.
For the reader's convenience, we recall the key definitions and results from \cite{australia} below. 

Let $\GG$ be an ultragraph. Define the set \[X=\mathfrak{p}^\infty\cup\{(\alpha,v): \alpha \in \mathcal{G}^*, |\alpha|\geq 1, v \in G^0_s \cap r(\alpha) \}\cup\{(v,v): v\in G^0_s \}.\]
We extend the range and source maps defined on $\GG$ to elements $(\alpha,v)\in X$ by defining $r(\alpha,v)=v$ and $s(\alpha,v)=s(\alpha)$. Furthermore, we extend the length map to the elements $(\alpha,v)$ by defining $|(\alpha,v)|:=|\alpha|$.

Let $\F$ denote the free group generated by $\mathcal{G}^1$, and $W\subseteq \F$ be the set of paths in $\GG^*$ with strictly positive length. Now, define the following subsets of $X$:
\begin{itemize}
\item for $a\in W$, let $X_a=\{x\in X:x_1..x_{|a|}=a\}$;

\item for $b\in W$, let $X_{b^{-1}}=\{x\in X:s(x)\in r(b)\}$;
\item for $a,b\in W$ with $r(a)\cap r(b)\neq \emptyset$, let $$X_{ab^{-1}}=\left \{x\in X:|x|>|a|, \,\,\,  x_1...x_{|a|}=a \text{ and }s(x_{|a|+1})\in r(b)\cap r(a)\right\}\bigcup$$ $$\bigcup\left\{(a,v) \in X:v\in r(a)\cap r(b)\right\};$$
\item for the neutral element $0$ of $\F$, let $X_0=X$;
\item for all the other elements $c$ of $\F$, let $X_c=\emptyset$.
\end{itemize}

Furthermore, for each $A\in \mathcal{G}^0$ and $b\in W$, let 
$$X_A=\{x\in X:s(x)\in A \}$$ and 
$$X_{bA}=\{x\in X_b:|x|>|b| \text{ and } s(x_{|b|+1})\in A\}\cup\{(b,v)\in X_b:v\in A\}.$$

Next, we define the usual erase and add maps:
for $a\in W$ define $\theta_a:X_{a^{-1}}\rightarrow X_a$ by $ \theta_a(x)= ax$, if $|x|=\infty$, and $ \theta_a((\alpha,v))= (a\alpha,v)$. Of course,  $\theta_a^{-1}:X_a\rightarrow X_{a^{-1}}$ is defined as the inverse of $\theta_a$; for $a,b\in W$ define $\theta_{ab^{-1}}:X_{ba^{-1}}\rightarrow X_{ab^{-1}}$ by $ \theta_{ab^{-1}}(x)= ay$ if  $|x|=\infty$, $ \theta_{ab^{-1}}((b\alpha,v))=
(a\alpha,v)$, and $ \theta_{ab^{-1}}((b,v))=
(a,v)$. For the neutral element $0\in \F$ define $\theta_0:X_0\rightarrow X_0$ as the identity map. Finally, for all the other elements $c$ of $\F$ define $\theta_c:X_{c^{-1}}\rightarrow X_c$ as the empty map.

The above maps together with the subsets $X_t$ form a partial action of $\F$ on $X$, that is $(\{\theta_t\}_{t\in \F}, \{X_t\}_{t\in \F})$ is such that $X_0=X$, $\theta_0=Id_x$, $\theta_c(X_{c^{-1}}\cap X_t)=X_{ct}\cap X_c$ and $\theta_c\circ\theta_t=\theta_{ct}$ in $X_{t^{-1}}\cap X_{t^{-1}c^{-1}}$. This partial action induces a partial action on the level of the $R$-algebra of functions (with point-wise sum and product) $F(X)$. 
More precisely, let $D$ be the subalgebra of $F(X)$ generated by all the finite sums of all the finite products of the characteristic maps $\{1_{X_A}\}_{A\in \mathcal{G}^0}$, $\{1_{bA}\}_{b\in W, A\in \mathcal{G}^0}$ and $\{1_{X_c}\}_{c\in \F}$. Also define, for each $t\in \F$, the ideals $D_t$ of $D$, as being all the finite sums of finite products of the characteristic maps $\{1_{X_t}1_{X_A}\}_{A\in \mathcal{G}^0}$, $\{1_{X_t}1_{bA}\}_{b\in W, A\in \mathcal{G}^0}$ and $\{1_{X_t}1_{X_c}\}_{c\in \F}$. Now, for each $c\in \F$, define the $R$-isomorphism $\beta_c:D_{c^{-1}}\rightarrow D_c$ by $\beta_c(f)=f\circ \theta_{c^{-1}}$. Then $(\{\beta_t\}_{t\in \F}, \{D_t\}_{t\in \F})$ is a partial action of $\F$ on $D$. 

\begin{remark}\label{spanDt} From now on we will use the notation $1_A$, $1_{bA}$ and $1_t$ instead of $1_{X_A}$, $1_{X_{bA}}$ and $1_{X_t}$, for $A\in \mathcal{G}^0$, $b\in W$ and $t\in \F$. Also, we have the following description of the ideals $D_t$:
$$D=D_0=\text{span}\{1_A, 1_c, 1_{bA}:A\in \mathcal{G}^0, c\in \F\setminus\{0\}, b\in W\}, $$ and, for each $t\in \F$, 
$$D_t=\text{span}\{1_t1_A, 1_t1_c, 1_t1_{bA}:A\in \mathcal{G}^0, c\in \F, b\in W\}. $$ 
\end{remark}

The partial skew group ring associated with the partial action defined above is isomorphic to the ultragraph Leavitt path algebra, as we recall below.

\begin{teorema}\cite[Theorem 2.8]{reduction}\label{isom} Let $\mathcal{G}$ be an ultragraph, $R$ be a unital commutative ring, and $L_R(\mathcal{G})$ be the Leavitt path algebra of $\mathcal{G}$. Then there exists an $R$-isomorphism $\phi:L_R(\mathcal{G})\rightarrow D\rtimes_\beta \F$ such that $\phi(p_A)=1_A\delta_0$, $\phi(s_e^*)=1_{{e^{-1}}}\delta_{e^{-1}}$ and $\phi(s_e)=1_e\delta_e$ for each $A\in \mathcal{G}^0$ and $e\in \mathcal{G}^1$.
\end{teorema}

Notice that $D\rtimes_\beta \F$ is $\F$-graded, with grading $\{D_t\delta_t\}_{t\in \F}$, and so, via the isomorphism of Theorem~\ref{isom}, we get an $\F$-grading on $L_R(\mathcal{G})$.

We now proceed to characterize which ultragraph Leavitt path algebras are strongly $\F$-graded and which are epsilon-strongly $\F$-graded. As we see below, asking for a strong $\F$-grading on an ultragraph Leavitt path algebra is very restrictive.

\begin{teorema}\label{voandonaagua}
Let $\mathcal{G}$ be an ultragraph. Then, $L_R(\mathcal{G})$ is strongly $\F$-graded if, and only if, $\mathcal{G}$ has only one edge $e$, in which case $s(e)=r(e)$. 
\end{teorema}

\demo Suppose that $L_R(\mathcal{G})\cong D \rtimes_\beta \F$ is strongly $\F$-graded. 

Let $e$ and $f$ be edges in $\mathcal{G}$.  Notice that each element of $D_e\delta_eD_{e^{-1}}\delta_{e^{-1}}$ is of the form $1_eg\delta_0$, with $g\in D_0$, and so $1_e\delta_0$ is a left unit for $D_e\delta_eD_{e^{-1}}\delta_{e^{-1}}$. Furthermore, $1_f\delta_0\in D\delta_0 = D_e\delta_eD_{e^{-1}}\delta_{e^{-1}}$ and hence $1_e \delta_0 1_f \delta_0 = 1_f \delta_0$. Since $1_f\delta_0\neq 0$, we obtain that $e=f$.

Denote by $e$ the unique edge in $\GG$. Similarly to what was done above, we obtain that 
$1_{r(e)}\delta_0$ is a left unit for $D_{e^{-1}}\delta_{e^{-1}}D_e\delta_e$.  Let $v$ be a vertex in $\GG^0$.  Then $1_v\delta_0\in D_0= D_{e^{-1}}\delta_{e^{-1}}D_e\delta_e$ and, since $1_v\delta_0 \neq 0$, we get that $v\in r(e)$. In particular, we obtain that $s(e)\in r(e)$ and so $ee$ is a finite path in $\GG$.

To finish this part of the proof it remains to show that $s(e)=r(e)$. Suppose that $|r(e)|>1$. Then, there exists a vertex $v\in r(e)$ which is a sink, and so $(e,v)\in X$. Notice, as before, that $1_{ee}\delta_0$ is a left unit for $D_{ee}\delta_{ee}D_{(ee)^{-1}}\delta_{(ee)^{-1}}=D_0\delta_0$. But $1_e \delta_0\in D_0 \delta_0$ and $1_e(e,v)\neq 0$, while $1_{ee}(e,v)=0$, leading to a contradiction.


For the converse, suppose that $\mathcal{G}$ has only one edge $e$ and $r(e)=s(e)$. Then $D_g=D_0$ for each $g\in \F$ (since $X$ contains only one point and, moreover, $X_g=X$ for each $g\in \F$). It is clear that  $D_g\delta_gD_h\delta_h\subseteq D_{gh}\delta_{gh}$. For the other inclusion, let  $g,h\in \F$ and $b\in D_{gh}$. Then $b\delta_{gh}=1_X\delta_gb\delta_h$, and so $D_{gh}\delta_{gh}\subseteq D_g\delta_gD_h\delta_h$. Therefore $D\rtimes_{\beta}\F$ is strongly $\F$-graded.

\fim

\begin{remark}
When dealing with $\F$-grading we are implicitly assuming that the ultragraph has at least one edge. Otherwise, the ultragraph consists of a collection of sinks and there is no grading over the free group of the edges.
\end{remark}

\begin{teorema}\label{ornitorrinco}
Let $\mathcal{G}$ be an ultragraph and let $\{D_t\delta_t\}_{t\in \F}$ as above. With this grading,  $L_R(\mathcal{G})$ is epsilon-strongly $\F$-graded if, and only if, $G^0\in \mathcal{G}^0$.
\end{teorema}

\demo If $L_R(\mathcal{G})$ is epsilon-strongly $\F$-graded then, by Proposition~\ref{G^0 unit}, $G^0\in \mathcal{G}^0$. 

For the converse, suppose that $G^0\in \mathcal{G}^0$. We will use the second item of Proposition~\ref{epsilon-strongly} to show that $D\rtimes_\beta \F$ is epsilon-strongly $\F$-graded (and hence, by Theorem~\ref{isom}, it follows that $L_R(\mathcal{G})$ is epsilon-strongly $\F$-graded).

We prove first that, for each $t\in \F$, there exists a unit for $D_t\delta_t D_{t^{-1}}\delta_{t^{-1}}$. For $t=0$, $1_{G^0}\delta_0$ is a unit of $D_0\delta_0D_0\delta_0$, since $1_{G^0}$ is a unit of $D_0$. Let $t\in \F\setminus \{0\}$ and $h\in D_t$. Then, $1_t\delta_0 h\delta_t=1_th\delta_t=h\delta_t$ (where the last equality holds since $1_t$ is a unit of $D_t$) and hence $1_t\delta_0=1_t\delta_t1_{t^{-1}}\delta_{t^{-1}}$ is a unit of $D_t\delta_tD_{t^{-1}}\delta_{t^{-1}}$. Similarly, one shows that for each $g\in D_{t^{-1}}$, $g\delta_{t^{-1}}1_t\delta_0=g 1_{t{^{-1}}} \delta_{t^{-1}}=g\delta_{t^{-1}}$.
It follows that $1_t\delta_0$ is a unit of  $D_t\delta_tD_t\delta_{t^{-1}}$. 

It remains to show that $D_t\delta_tD_{t^{-1}}\delta_{t^{-1}}D_t\delta_t=D_t\delta_t$, for all $t\in \F$. It is clear that $D_t\delta_tD_{t^{-1}}\delta_{t^{-1}}D_t\delta_t\subseteq D_t\delta_t$. For the other inclusion, recall that $1_t\delta_0\in D_t\delta_t D_{t^{-1}}\delta_{t^{-1}}$. Then, for $h\in D_t$, we have that $$h\delta_t=1_t\delta_0 h\delta_t\in D_t\delta_tD_{t^{-1}}\delta_{t^{-1}}D_t\delta_t,$$ and therefore $D_t\delta_t\subseteq D_t\delta_tD_{t^{-1}}\delta_{t^{-1}}D_t\delta_t$ as desired.

\fim 

It is proved in \cite[Lemma~6.11]{Imanfar} that an ultragraph Leavitt path algebra $L_R(\GG)$ is unital if, and only if, $G^0\in \GG^0$. We, therefore, have the following characterization of epsilon-strongly $\F$-graded ultragraph Leavitt path algebras.

\begin{corolario}\label{platapus}
Let $\GG$ be an ultragraph. Then the associated Leavitt path algebra is epsilon-strongly $\F$-graded if, and only if, it is unital.
\end{corolario}
 
\demo
This follows Theorem~\ref{ornitorrinco} and \cite[Lemma~6.11]{Imanfar}.

\fim

\vspace{1.5pc}

Daniel Gon\c{c}alves, Departamento de Matem\'{a}tica, Universidade Federal de Santa Catarina, Florian\'{o}polis, 88040-900, Brazil.

Email: daemig@gmail.com

\vspace{0.5pc}
Danilo Royer, Departamento de Matem\'{a}tica, Universidade Federal de Santa Catarina, Florian\'{o}polis, 88040-900, Brazil.

Email: daniloroyer@gmail.com
\vspace{0.5pc}

\end{document}